\documentclass[10pt]{amsart}
\usepackage[psamsfonts]{amssymb}
\usepackage{amsthm,amsmath}
\usepackage{color}

\title[Invariant functionals on completely distributive lattices]{Invariant functionals on completely\\ distributive lattices}

\author{Marta Cardin}
\address{Department of Applied Mathematics, University Ca' Foscari of Venice \\
Dorsoduro 3825/E--30123, Venice, Italy}
\email{mcardin[at]unive.it }

\author{Miguel Couceiro}
\address{Mathematics Research Unit, University of Luxembourg \\
6, rue Richard Coudenhove-Kalergi,
L-1359 Luxembourg, Luxembourg }
\email{miguel.couceiro[at]uni.lu }

\date{February 2, 2010}

\begin{document}

\theoremstyle{plain}
\newtheorem{theorem}{Theorem}
\newtheorem{lemma}[theorem]{Lemma}
\newtheorem{proposition}[theorem]{Proposition}
\newtheorem{corollary}[theorem]{Corollary}
\newtheorem{fact}[theorem]{Fact}
\newtheorem{claim}{Claim}
\newtheorem*{main}{Main Theorem}

\theoremstyle{definition}
\newtheorem{definition}[theorem]{Definition}
\newtheorem{example}[theorem]{Example}

\theoremstyle{remark}
\newtheorem{remark}{Remark}

\newcommand{\N}{\mathbb{N}}                     % pos. integers
\newcommand{\R}{\mathbb{R}}                     % reals
\newcommand{\co}[1]{\ensuremath{\overline{#1}}}
\newcommand{\vect}[1]{\ensuremath{\mathbf{#1}}} % vector
\def\med{\mathop{\rm med}\nolimits}

\begin{abstract} In this paper we are interested in functionals defined on completely distributive lattices and
which are invariant under mappings preserving {arbitrary} joins and meets.
We prove that the class of nondecreasing invariant functionals coincides with the class of Sugeno integrals
associated with $\{0,1\}$-valued capacities, the so-called term functionals, thus extending previous results both 
to the infinitary case as well as to the realm of completely distributive lattices. Furthermore, we show that, in the case of functionals 
over complete chains, the nondecreasing condition is redundant.
Characterizations of the class of Sugeno integrals, as well as its superclass comprising all polynomial functionals,
are provided by showing that the axiomatizations (given in terms of homogeneity) of their restriction to finitary functionals still hold  
over completely distributive lattices. 
We also present canonical normal form representations of polynomial functionals on completely distributive lattices,
which appear as the natural extensions to their finitary counterparts, and as a by-product we obtain an axiomatization
of complete distributivity in the case of bounded lattices.
\end{abstract}


\maketitle

%---------------------------------------------------------------------------------------------- Section 1
\section{Introduction}
The process of merging or combining sets of values (often numerical) into a single one is usually achieved by the
so-called aggregation functionals. A classical example of such an aggregation functional is that of a weighted arithmetic mean.
The importance of aggregation functionals is made apparent by their wide use, not only in pure mathematics
(e.g. in the theory of functional equations, measure and integration theory), but also in several applied fields such as  operations research, computer and information sciences, economics and social sciences, as well as in other experimental areas of physics and natural sciences.
The growing need to fuse several inputs into a single output in such a way that the resulting value somehow represents all the inputs led to the theory of aggregation whose main problem is to propose and describe aggregation functionals suitable for a required application. For general background, see \cite{BelPraCal07, GraMurSug00} and for a recent reference, see \cite{GraMarMesPap09}.

There are two main approaches to describing classes of aggregation functionals. One defines a class by explicitly specifying its members (as in the case of arithmetic means or certain fuzzy integrals such as Sugeno integrals), the other proposes properties and conditions which are fulfilled by exactly the members of the class, that is, which axiomatize the given class of aggregation functionals. Both are equally important and interesting in their own right.
In the former approach, such descriptions are many times given in terms of normal form representations (e.g., polynomial or disjunctive normal forms). In the latter approach, such characterizations are given in terms of necessary and sufficient conditions often expressed in the form of functional equations.

Let ${\bf X}$ be a set, $A$ be a nonempty set, and $K$ a set of mappings $\varphi \colon {\bf X} \to {\bf X}$.
In this paper we are particularly interested in classes of aggregation functionals $F\colon {\bf X}^A\to {\bf X}$ which satisfy functional equations of the following type
\begin{equation}\label{eq:1}
F(\varphi\circ f)= \varphi\circ F(f),
\end{equation}
where $f$ ranges over the elements of ${\bf X}^A$ and $\varphi$ over the members of $ K$. Functionals fulfilling (\ref{eq:1}) are said to be invariant under $K$ (or simply, $K$-invariant).  In the framework of utility theory,  (\ref{eq:1}) expresses the fact that aggregating the transformed inputs is the same as transforming the aggregated inputs.

In the case when ${\bf X}$ is a dense and linearly ordered set, or dense chain (in fact, a real interval), 
$A$ is finite  and $K$ is the set of increasing bijections (automorphisms) on ${\bf X}$, solutions of (\ref{eq:1}) were obtained under certain conditions by Greco \cite{Greco}, Ovchinnikov \cite{OV}, Marichal \cite{Mar} and Marichal and Mathonet \cite{MarMath}, and shown to coincide with certain (discrete) Sugeno integrals determined by $\{0,1\}$-valued capacities. These results were then extended by Dukhovny and Ovchinnikov \cite{OvDuck02} to the case when $A$ is infinite.  In a somewhat different setting, {Chambers} \cite{Cham} considers invariance under continuous and strictly increasing mappings on real numbers. 
%if you want to go beyond this reference you should provide sentences that make sense and give an idea to the reader. Not something that just confuses the reader.
 For further variants and developments see \cite{OvDuck, Mar02, MesRuck}, and for a recent survey on the topic, see Marichal and Mesiar \cite{MM}.

In this paper we seek the solutions of (\ref{eq:1}) in a rather different setting. Letting $A$ be an arbitrary nonempty set (and thus subsuming infinitary mappings), we study functionals $F\colon L^A\to L$ defined and valued over a completely distributive lattice $L$ (not necessarily dense nor linearly ordered).  Of particular interest are those functionals
defined by fuzzy integrals with respect to certain classes of capacities (or fuzzy measures, in the terminology of \cite{Sug74}). Also, instead of automorphisms, we consider invariance under mappings which preserve arbitrary meets and joins. 
We show that, under nondecreasing monotonicity, invariant functionals coincide exactly with (infinitary) Sugeno integrals associated with $\{0,1\}$-valued capacities, i.e., idempotent aggregation functionals which can be obtained from projections and constants by making use of arbitrary meets and joins. By relaxing the invariance conditions, we obtain characterizations for the class of Sugeno integrals as well as for the superclass of lattice polynomial functionals.
To this extent, we provide canonical representations
of polynomial functionals on completely distributive lattices, and present a characterization of complete distributivity, alternative to those given
by Raney \cite{Ra} and Tunnicliffe \cite{Tunn}.
Moreover, we prove that
in the case when $L$ is a complete chain (not necessarily dense), nondecreasing monotonicity follows from invariance.

The structure of the paper is as follows. We start by recalling some well-known concepts and terminology in lattice theory.
In Section 3, we focus on functionals over completely distributive lattices. We introduce the notion of complete distributivity as given in \cite{Ra} and present a characterization of this property due to Tunnicliffe \cite{Tunn}. As it will become clear (Subsection 3.2), this stronger distributive law
 enables neat normal form representations of polynomial functionals, in particular, of Sugeno integrals and term functionals (i.e., Sugeno integrals associated with $\{0,1\}$-valued capacities). Using canonical  representations of the latter  we obtain an alternative characterization
 of complete distributivity  (Subsection 3.3). In Section 4, we study nondecreasing functionals invariant under mappings preserving arbitrary meets and joins (which we refer to as continuous) and show, in Subsection 4.2, that they coincide exactly with term functionals (Sugeno integrals  associated with $\{0,1\}$-valued capacities). Moreover, in the particular case when $L$ is a complete chain, we prove that nondecreasing monotonicity is a redundant condition since it follows from invariance.
 The superclasses of Sugeno integrals and, more generally, of polynomial functionals, are considered in Subsection 4.3 and shown to be axiomatized by weaker forms of invariance, namely, by homogeneity.

\section{Basic notions and terminology}
A \emph{lattice} is an algebraic structure $\langle L,\wedge, \vee \rangle$ where $L$ is a nonempty set, called \emph{universe}, and where the
two binary operations $\wedge$ and $\vee$ satisfy the commutative, the associative, the absorption, and the idempotent laws. With no danger of
ambiguity, we will denote lattices by their universes. A lattice $L$ is said to be \emph{distributive} if, for every $a,b,c\in L$,
$$
a\vee (b\wedge c)= (a\vee b)\wedge (a\vee c)\quad \mbox{or, equivalently,} \quad a\wedge (b\vee c)= (a\wedge b)\vee (a\wedge c).
$$

For $x,y\in L$, $x\leqslant y$ simply means that $x\wedge y=x$ or, equivalently, $x\vee y=y$.
A subset $S$ of a lattice $L$ is said to be \emph{convex} if for every $a,b\in S$ and every $c\in L$ such that $a\leqslant
c\leqslant b$, we have $c\in S$. For any subset $S\subseteq L$, we denote by $\co{S}$ the convex hull of $S$, that is, the smallest convex
subset of $L$ containing $S$.
A \emph{chain} is a lattice such that for every $a,b\in L$ we have  $a\leqslant b$ or $b\leqslant a$. Clearly, every chain is distributive.

For an arbitrary nonempty set $A$ and a lattice $L$, the set $L^A$ of all functions from $A$ to $L$ also constitutes a lattice under the
operations
$$
(f\wedge g)(x)=f(x)\wedge g(x) \quad \textrm{ and } \quad (f\vee g)(x)=f(x)\vee g(x),
$$
for every $f,g\in L^A$. We denote the elements of $L$ by lower case letters $a,b,c,\ldots$, and the elements of $L^A$ by usual function symbols
$f,g,h, \ldots$.
A lattice $L$ is said to be \emph{bounded} if it
has a least and a greatest element, usually denoted by $0$ and $1$, respectively. We use $\vect{0}$ and $\vect{1}$ to denote the
least element and greatest element, respectively, of $L^A$. Likewise and with no danger of ambiguity, for each $c\in L$, we denote by $\vect{c}$ the constant c map in $L^A$. If $L$ is bounded, then for each $X\subseteq A$, we denote by $I_X$ the \emph{characteristic function} of $X$ in $L^A$, i.e.,
\[
I_X(x) =
\begin{cases}
1, & \text{if $x\in X$} \\
0, & \text{otherwise.}
\end{cases}
\]
For further background in lattice theory we refer the
reader to, e.g., Birkhoff \cite{Bir}, Davey and Priestley~\cite{DP}, Gr\"atzer~\cite{Grae03}, and Rudeanu~\cite{Rud01}.

By a (\emph{lattice}) \emph{functional} on $L$ we mean a mapping $F\colon L^A\to L$, where $A$ is a nonempty set.
The \emph{range} of a functional  $F\colon L^A\to L$ is defined by $\mathcal{R}_F=\{F(f) : f\in L^A\}$. A functional
 $F\colon L^A\to L$ is said to be \emph{nondecreasing} if, for every $f, g\in L^A$ such that
$f(i)\leqslant g(i)$, for every $i\in A$, we have $F(f)\leqslant F(g)$.
Note that if $F$ is nondecreasing, then $\co{\mathcal{R}}_F=[F(\vect{0}),F(\vect{1})]$.
An \emph{aggregation functional} on a bounded lattice $L$ is a nondecreasing functional $F\colon L^A\to L$ such that $\co{\mathcal{R}}_F=L$, that is, $F(\vect{c})=c$ for $c\in \{0,1\}$. For instance, for each $a\in A$, the projection $F_a\colon L^A\to L$ defined by
$F_a(f)=f(a)$, is an aggregation functional.

As mentioned, in this paper we are particularly interested in certain lattice functionals (in particular, aggregation functionals) which bring the notion of lattice polynomial function (see e.g. \cite{Goo67,CouMar})
to the infinitary function setting. To this extent we need to consider some special bounded lattices.
A lattice $L$ is said to be \emph{complete} if $\bigwedge S=\bigwedge_{x\in S}x$ and $\bigvee S=\bigvee_{x\in S}x$ exist for every $S\subseteq L$. Clearly, every complete lattice is also bounded.

Let $L$ be a complete lattice and $A$ a nonempty set. By a (\emph{lattice}) \emph{polynomial functional} on $L$ we mean a functional $F\colon L^A\to L $ which can be obtained by the following rules:
\begin{enumerate}
\item[(i)] Each projection $F_a$, $a\in A$, and each  constant $c$ functional, $c\in L$, is a polynomial functional;
\item[(ii)] If $\mathcal{F}=\{F_i\colon L^A\to L:i\in I\}$ is a family of polynomial functionals, then $\bigwedge \mathcal{F}$ and $\bigvee \mathcal{F}$ are polynomial functionals.
\end{enumerate}
If $A$ is finite, then polynomial functionals are usually referred to as (lattice) polynomial functions.
Polynomial functionals $F\colon L^A\to L $ which are \emph{idempotent}, that is, satisfying
$$F(\vect{c})=c, \quad \mbox{for every constant map $\vect{c}\in L^A$,}$$
are called \emph{Sugeno integrals} on $L$.
(Sugeno integrals were introduced by Sugeno \cite{Sug74, Sug77} on linearly ordered domains.
 In the finitary case, Marichal \cite{Mar08} observed that this concept can be extended to the setting of bounded distributive lattices
 by defining Sugeno integrals as idempotent polynomial functions.)
We shall refer to polynomial functionals obtained from projections and rule (ii) as \emph{term functionals}. In other words, term functionals are exactly those polynomial functionals satisfying $F(f)\in\{0,1\}$, for every $f\in \{0,1\}^A$, i.e., Sugeno integrals associated with $\{0,1\}$-valued capacities.

{\flushleft \bf Fact.}
\emph{Every polynomial functional is nondecreasing and thus every Sugeno integral is an aggregation functional.}

\section{Polynomial functionals on completely distributive lattices}

In this section, we are interested in polynomial functionals on complete lattices which satisfy a stronger variant of  distributivity, namely, complete distributivity. As we will see, this more stringent form of distributivity allows neat representations of polynomial functionals.
We start by recalling the notion of complete distributivity and present a characterization of this property due to Tunnicliffe \cite{Tunn}.
 Then we provide normal form representations of polynomial functionals, in particular, of Sugeno integrals and term polynomials, on completely distributive lattices. As it will become clear (Subsection 3.3), the latter leads to an alternative characterization of complete distributive lattices given with respect to certain term functionals.

\subsection{General background on completely distributive lattices}

A complete lattice $L$ is said to be \emph{completely distributive} if the following more stringent distributive law holds
\begin{equation}\label{Eq:distrib}
\bigwedge_{i\in I}\Bigl(\bigvee_{j\in J} x_{ij}\Bigr)=\bigvee_{f\in J^{I}}\Bigl(\bigwedge_{i\in I} x_{if(i)}\Bigr),
\end{equation}
for every doubly indexed subset $\{x_{ij}: i\in I, \, j\in J\}$ of $L$. Note that every complete chain is completely distributive.
Complete distributivity is a self-dual property. This was observed by Raney \cite{Ra}
 who showed that (\ref{Eq:distrib}) and its dual are equivalent, and thus that either is sufficient to define complete distributivity.

 In this paper we shall make use of an alternative characterization of complete distributivity due to Tunnicliffe \cite{Tunn},
 which essentially relies on the notion of ``cross-cuts". Let $\mathcal H$ and $ \mathcal K$ be two nonempty families of subsets of
 a complete lattice $L$. A pair $(\mathcal H, \mathcal K)$ is called a \emph{cone}, if $H\cap K\neq \emptyset$
 for all $H \in \mathcal H, \,K \in \mathcal K$. Cones can be naturally ordered with respect to inclusion. More precisely, this order is given by
 $ (\mathcal H, \mathcal K)\leqslant (\mathcal H^{\prime}, \mathcal K^{\prime})$  if $\mathcal H
\subseteq \mathcal H^{\prime} $ and $\mathcal K \subseteq \mathcal K^{\prime}$.  A \emph{ultracone} is simply
% here "simply" is just to improve readibility since it sounds funny "A ultracone is a cone". But again this is a matter of taste we can delete it
a cone $(\mathcal H, \mathcal K)$
which is maximal with respect to this ordering, i.e., if $ (\mathcal H, \mathcal K)\leqslant (\mathcal H^{\prime}, \mathcal K^{\prime})$, then
$\mathcal H= \mathcal H^{\prime} $ and $\mathcal K = \mathcal K^{\prime}$.

\begin{theorem}\label{Th:Tunn}\emph{(\cite{Tunn})}
Let $L$ be a complete lattice. Then  $L$ is completely distributive if
and only if
$$\bigvee_{H \in \mathcal H}\bigwedge H=\bigwedge_{K \in \mathcal K}\bigvee K,$$
for every ultracone $(\mathcal H, \mathcal K)$.
\end{theorem}

\subsection{Representations of polynomial functionals on completely distributive lattices}

As mentioned, when considered on completely distributive lattices, lattice polynomial functionals have neat normal form representations.
The following result provides a canonical construction of such representations. 
%extending to the infinitary case the normal form representations of polynomial functions obtained by Goodstein \cite{Goo67}.

\begin{theorem} \label{Prop:2}
Let $L$ be a completely distributive lattice, $A$ an arbitrary
nonempty set, and let $F \colon  L^{A}\rightarrow L $ be a nondecreasing functional.
Then, for every $f\in L^A$,
$$P_F(f)= \bigvee_{X\in\mathcal A} F(I_X) \wedge \bigwedge_{x\in X} f(x)=\bigwedge_{X\in\mathcal B} F(I_{A\setminus X}) \vee \bigvee_{x\in X} f(x)=P^F(f),$$
where $\mathcal A=\{X\subseteq A\colon F(I_X)\neq 0\}$ and $\mathcal B=\{X\subseteq A\colon F(I_{A\setminus X})\neq 1\}$. In particular, if
$F$ is a polynomial functional, then $F=P_F=P^F$.
\end{theorem}

\begin{proof} Let $\mathcal{H}=\{H\subseteq L: H\supseteq f(X)\cup\{F(I_Y):X\subseteq Y\}  \,\mbox{ for some }  \, X\subseteq A\}$ and
$\mathcal{T}=\{K\subseteq L: K\supseteq f(X)\cup\{F(I_{A\setminus X})\}  \,\mbox{ for some }  \, X\subseteq A\}$.
Observe that
$$\bigvee_{H \in \mathcal H}\bigwedge H=P_F(f)\quad \mbox{ and } \quad \bigwedge_{K \in \mathcal{T}}\bigvee K =P^F(f).$$

 Let $\mathcal{K}= \{K\subseteq L : K\cap H\neq \emptyset \,\mbox{ for every }\, H\in \mathcal{H} \}$.
 We claim that $(\mathcal{H}, \mathcal{K})$ is a ultracone. Suppose that $(\mathcal{H}', \mathcal{K}')$ is a cone such that
 $(\mathcal{H}, \mathcal{K})\leqslant (\mathcal{H}', \mathcal{K}')$.
 For every $X\subseteq A$, we have $f(X)\cup\{F(I_Y):X\subseteq Y\} \in \mathcal{H}\subseteq \mathcal{H}'$. Hence, if $K\in \mathcal{K}'$, then
 $K\cap (f(X)\cup\{F(I_Y):X\subseteq Y\})\neq \emptyset  $ and $K\in \mathcal{K}$.

 Now, for the sake of a contradiction, suppose that  there is $H\in \mathcal{H}'\setminus \mathcal{H}$. Then, for every $X\subseteq A$, there is $Z_X\subseteq X $ such that  $f(Z_X)\not \subseteq H$ or there is $W_X\supseteq X$ such that $F(I_{W_X})\not \in H$.
 Let $B=\{f(z)\not\in H:z\in Z_X, \, X\subseteq A\}\cup\{F(I_{W_X}): X\subseteq A\}$. Then $B\in \mathcal{K}=\mathcal{K}'$, but $B\cap H=\emptyset$. This yields the desired contradiction.

Now, since  $(\mathcal{H}, \mathcal{K})$ is a ultracone, it follows from Theorem \ref{Th:Tunn} that
$$\bigwedge_{K \in \mathcal{K}}\bigvee K=P_F(f).$$
 Thus to complete the proof, it is enough to show that
  $$\bigwedge_{K \in \mathcal{K}}\bigvee K\leqslant P^F(f)\quad \mbox{ and }\quad P^F(f)\leqslant P_F(f).$$

 Let $K\in \mathcal{T}$ and fix $X\subseteq A$ such that $K\supseteq f(X)\cup\{F(I_{A\setminus X}\}$.
 If for every $Y\subseteq A$, $K\cap f(Y)\neq \emptyset$, then we have $K\in \mathcal{K}$.
 If there is $Y\subseteq A$ such that $K\cap f(Y)=\emptyset$, then $f(X)\cap f(Y)=\emptyset$, and hence $X\cap Y=\emptyset$. Therefore, $X\subseteq A\setminus Y$ and thus
 $K \cap \big( f(Y)\cup\{F(I_Z):Y\subseteq Z\}\big)\neq \emptyset$. Since this holds for each such $Y$, we have that $K\in \mathcal{K}$ and so
 $\bigwedge_{K \in \mathcal{K}}\bigvee K\leqslant P^F(f)$.

 Since $L$ is completely distributive, we can find $G\colon L^A\to L$ such that
 $$\bigvee_{X\subseteq A} G(I_X) \wedge \bigwedge_{x\in X} f(x)= P^F(f).$$
 Indeed, we have $$P^F(f)= \underset{X\in\mathcal B}{\bigwedge}F(I_{A\setminus X}) \vee 
 \underset{x\in X}{\bigvee} f(x)=  \underset{X\in\mathcal{P} (A)}{\bigwedge}
 \underset{x\in A}{\bigvee} a_{x,X},$$
 where $\mathcal{P} (A)$ denotes the power set of $A$, and $a_{x,X}=F(I_{A\setminus X})\vee f(x)$ if $x\in X$, and $a_{x,X}=0$ otherwise.
 By (\ref{Eq:distrib}), we then have
 \begin{eqnarray*}
 P^F(f)&=& \bigwedge_{X\in\mathcal{P} (A)} \, \bigvee_{x\in A} a_{x,X} \, = \, 
 \bigvee_{f\in A^{\mathcal{P} (A)}} \, \bigwedge_{X\in\mathcal{P} (A)}a_{f(X),X}\\
  &=& \bigvee_{f\in A^{\mathcal{P} (A)}} \, \bigwedge_{y\in Im(f)}\, \,  \bigwedge_{Y\colon y=f(Y)}a_{y,Y}\\
 &=& \bigvee_{f\in A^{\mathcal{P} (A)}} \, \bigwedge_{y\in Im(f)}\,  \, 
 \underset{y\in Y}{\bigwedge_{Y\colon y=f(Y)}}a_{y,Y}\\
 % &=& \bigvee_{f\in A^{\mathcal{P} (A)}}\bigwedge_{y\in Im(f)}\bigwedge_{y=f(Y), y\in Y}a_{y,Y}\\
  &=& \bigvee_{f\in A^{\mathcal{P} (A)}} \, \bigwedge_{y\in Im(f)}\,  \, 
  \underset{y\in Y}{\bigwedge_{Y\colon y=f(Y)}} \big(F(I_{A\setminus Y})\vee f(y)\big)\\
   &=& \bigvee_{f\in A^{\mathcal{P} (A)}} \, \bigwedge_{y\in Im(f)} (t_{y,1}\vee t_{y,2})
 \end{eqnarray*}
where $t_{y,1}=  \underset{y\in Y}{\underset{Y\colon y=f(Y)}{\bigwedge}}F(I_{A\setminus Y})$ and $t_{y,2}=f(y)$. Again by (\ref{Eq:distrib}),
we get 
$$
 P^F(f)   = \bigvee_{f\in A^{\mathcal{P} (A)}} \, \bigvee_{g\in \{1,2\}^{Im(f)}} \, \bigwedge_{y\in Im(f)} t_{y,g(y)}
 = \bigvee_{X\subseteq A} a_X \wedge \bigwedge_{x\in X} f(x),
$$
for a suitable choice of the $a_X$'s. This proves the claim, for we can define $G\colon L^A\to L$ by setting $G(f)=\underset{X\subseteq A}{\bigvee} a_X \wedge \underset{x\in X}{\bigwedge} f(x)$ for every $f\in L^A$.

Now, for every $X\subseteq A$,  we have $F(I_X)= P_F(I_X)=P^F(I_X),$
and so $$\bigvee_{Y\subseteq X}G(I_Y)=F(I_X).$$
 Since $F(I_X)=0$ for all $X\not\in\mathcal A $, we have that for every $f\in L^A$,
$$P^F(f)=\bigvee_{X\in\mathcal A} G(I_X) \wedge \bigwedge_{x\in X} f(x)\leqslant P_F(f).$$

For the last claim, we observe that $P_F(f)\leqslant F(f)\leqslant P^F(f),$ for all $f\in L^A$.
To see this let $X\subseteq A$, and set $c_X=\underset{x\in X}{\bigwedge} f(x)$ and $d_X=\underset{x\in X}{\bigvee} f(x)$.
If $F \colon  L^{A}\rightarrow L $ is a polynomial functional, then by complete distributivity we have 
$$
F(I_X)\wedge c_X\leqslant F(I_X\wedge \vect{c}_X)\leqslant F(f) \leqslant F(I_{A\setminus X}\vee \vect{d}_X) \leqslant F(I_{A\setminus X})\vee d_X.
$$ 
Since the above holds for every $X\subseteq A$, we get $F(f)=P_F(f)=P^F(f),$ for all $f\in L^A$.
Hence, the proof of Theorem \ref{Prop:2} is now complete.
 \end{proof}

\begin{remark}
\begin{itemize}
 \item[\,(1)] A functional $F \colon  L^{A}\rightarrow L $ is said to have a \emph{disjunctive normal form} (for short, DNF) if there exist $a_X\in L$, for every $X\subseteq A$, such that for every $f\in L^A$
\begin{equation}\label{DNF}
F(f)=\bigvee_{X\in\mathcal A} a_X\wedge \bigwedge f(X).
\end{equation}
Thus, from Theorem \ref{Prop:2} it follows that lattice polynomial functionals coincide exactly with those functionals  $F \colon  L^{A}\rightarrow L $ which have a DNF. 

\item[(2)] In the case when $A$ is finite, the same result holds with the condition of complete distributivity on $L$, relaxed to distributivity. This was first shown by Goodstein \cite{Goo67} where it was assumed that $L$ was both bounded and distributive. However, the boundness condition on $L$ is not really necessary, since $L$ can be extended to a bounded lattice $L'$ by adjoining a greatest element  $1$ and least element $0$,  if necessary, and we have: 
\begin{itemize}
\item[(i)] if $F$ is a polynomial functional over $L$ induced by a DNF as in (\ref{DNF}), then the same expression induces a polynomial functional $F'$ on $L'$ such that the restriction of $F'$ to $L$ coincides with $F$, and
\item[(ii)] if $F'$ is a polynomial functional (not constant $0$ nor $1$) on $L'$ represented as in (\ref{DNF}), then
 by omitting each term $a_X \wedge \bigwedge f(X)$ where $a_X = 0$, and replacing each term $a_X \wedge \bigwedge f(X)$ where $a_X = 1$ by $\bigwedge f(X)$, we obtain an equivalent DNF for $F'$, and the functional $F\colon L^A\to L$ induced by this new DNF coincides with the restriction of $F'$ to $L$. 
\end{itemize}
\item[(3)] Similarly, for an arbitrary set $A$, Theorem \ref{Prop:2} still holds when $L$ is a chain and the functionals considered are restricted to
$$S=\{f\in L^A\colon \mbox{ for }X\subseteq A, \, \bigwedge f(X) \mbox{ and } \, \bigvee f(X)\, \mbox{ exist}\}.$$
\end{itemize}
\end{remark}

Observe also that from Theorem \ref{Prop:2} it follows that every nondecreasing mapping $F\colon  \{0,1\}^A\to L$ can be extended to a polynomial functional on $L$, and that this extension is unique.
In the case of term functionals, Theorem \ref{Prop:2} reduces to the following result.

\begin{corollary}\label{cor:1}
  Let $L$ be a completely distributive lattice, $A$ an arbitrary
nonempty set, and let $F \colon  L^{A}\rightarrow L $ be a nondecreasing functional such that $F(I_X)\in \{0,1\}$ for every $X\subseteq A$.
Then, for every $f\in L^A$,
$$F(f)=P_\mathcal{A}(f)= \bigvee_{X\in\mathcal A} \bigwedge_{x\in X} f(x)=\bigwedge_{X\in\mathcal B} \bigvee_{x\in X} f(x)=P^\mathcal{B}(f),$$
where $\mathcal A=\{X\subseteq A\colon F(I_X)=1\}$ and $\mathcal B=\{X\subseteq A\colon F(I_{A\setminus X})=0\}$. 
\end{corollary}

\subsection{Alternative characterization of completely distributive lattices}

The description of term functionals given in Corollary \ref{cor:1} together with Theorem \ref{Th:Tunn},
lead to the following noteworthy characterization of complete distributivity.

\begin{theorem}\label{Thm:3.2}
A complete lattice $L$ is completely distributive  if and only if
 for every set $A$ and every
 family $\mathcal A$ of nonempty  subsets of $A$, we have
$$ P_{\mathcal A}(f)= P^{\mathcal B}(f) $$
 for every $f\in L^{A}$ where $\mathcal B= \{B\subseteq A : B\cap X\neq
\emptyset \,\mbox{ for all}\,\, X\in \mathcal A \}$.
\end{theorem}

\begin{proof}
Let $\mathcal A$ be a nonempty family of subsets of $A$,
 and set $\mathcal B= \{B\subseteq A : B\cap X\neq
\emptyset \,\mbox{ for all}\,\, X\in \mathcal A \}$.
Define $F \colon  L^{A}\rightarrow L $ by  $F(I_X)=1$, if there is $Y\in \mathcal{A}$ such that $Y\subseteq X$, and $F(I_X)=0$, otherwise.
Observe that $F$ is nondecreasing.

By Corollary \ref{cor:1}, to prove necessity it is enough to show that $\mathcal B= \{X\subseteq A : F(I_{A\setminus X})=0 \}$.
So suppose that $X\in \mathcal{B}$. Then, for every $Y\subseteq A\setminus X$, $Y\not \in \mathcal{A}$ and thus $F(I_{A\setminus X})=0$.
On the other hand, if $X\not\in \mathcal{B}$, then there is $Y\in \mathcal{A}$ such that $X\cap Y=\emptyset$.
Hence, $Y\subseteq A\setminus X$. Since $F$ is nondecreasing, we have $F(I_{A\setminus X})=1$, and thus
$X\not\in \{X\subseteq A : F(I_{A\setminus X})=0 \}$.

By Theorem \ref{Th:Tunn}, to show that the condition is sufficient it is enough to prove that for every ultracone $(\mathcal H, \mathcal K)$, we have
$$\bigvee_{H \in \mathcal H}\bigwedge H=\bigwedge_{K \in \mathcal K}\bigvee K.$$

So let  $(\mathcal H, \mathcal K)$ be a ultracone, and set
$$A= \{x\in L : x\in H \,\mbox{for some}\, H\in \mathcal H\}\cup \{x\in L : x\in K \,\mbox{for some}\, K\in \mathcal K\}.$$
Then, by setting $\mathcal{A}=\mathcal H$ and $\mathcal{B}=\mathcal K$ and taking  $f\in L^A$ as the identity $f(x)=x$, we have
$$ \bigvee_{H \in \mathcal H}\bigwedge H = P_{\mathcal A}(f)= P^{\mathcal B}(f)= \bigwedge_{K \in \mathcal K}\bigvee K.$$
 \end{proof}

\section{Invariant lattice functionals}

In this section, we consider a notion of invariance with respect to certain lattice homomorphisms which appear as the natural analogues of continuous maps on topological spaces.  As it turns out, it constitutes a suitable condition for axiomatizing classes of polynomial functionals.
We start by introducing this notion in Subsection 4.1 together with few observations relating it to homogeneity as introduced in \cite{FodRou95} and defined more generally in \cite{CouMar}. Then we provide characterizations of term functionals (Sugeno integrals associated with $\{0,1\}$-valued capacities) in the general setting of completely distributive lattices. This characterization is refined in the particular case of bounded chains 
by showing that nondecreasing monotonicity follows from invariance. 
In the last subsection,
 we focus on the superclasses of Sugeno integrals and polynomial functionals, and provide characterizations for these classes by 
relaxing invariance accordingly.

\subsection{Continuity and invariance}

Traditionally, continuous functions are defined as mappings which preserve certain limits.
Working on ordered structures whose operations reduce to meets and joins, the natural approach to continuity is to define it with respect to preservation of arbitrary meets and joins. As it turns out (see \cite{DP}, exercise 8.8), under a suitable topologization of complete lattices, the classical notion of continuity and that given below are equivalent.

Let $L$ be a complete lattice. A mapping $\gamma \colon L\to L$ is said to be \emph{continuous} if it preserves arbitrary meets and joins, i.e.,
for every $S\subseteq L$,
$$\gamma (\bigwedge S) =\bigwedge \gamma(S)\quad \mbox{and} \quad \gamma (\bigvee S) =\bigvee \gamma(S).$$
We  say that a functional $F \colon  L^{A}\rightarrow L $ on a completely distributive lattice $L$ is \emph{invariant} if, for every $f\in L^A$ and every continuous mapping $\gamma \colon L\to L$, we have
$$F(\gamma \circ f)=\gamma \circ F(f).$$
For instance, the \emph{median}
$\med(x_1,x_2,x_3)=(x_1\wedge x_2)\vee (x_1\wedge x_3)\vee (x_2\wedge x_3),$
is an invariant functional, since for every continuous mapping $\gamma \colon L\to L$, we have
\begin{eqnarray*}
\med\big(\gamma(x_1),\gamma(x_2),\gamma(x_3)\big)
&=&
(\gamma(x_1)\wedge\gamma(x_2))\vee (\gamma(x_1)\wedge\gamma(x_3))\\ 
&\vee& (\gamma(x_2)\wedge \gamma(x_3))\\
&=& \gamma(x_1\wedge x_2)\vee \gamma(x_1\wedge x_3)\vee \gamma(x_2\wedge x_3)\\
&=& \gamma\big((x_1\wedge x_2)\vee (x_1\wedge x_3)\vee (x_2\wedge x_3)\big)\\
&=&\gamma \circ \med(x_1,x_2,x_3).
\end{eqnarray*}
Observe that invariance subsumes the notions of homogeneity as considered in \cite{FodRou95, CouMar}. Indeed, a functional $F \colon  L^{A}\rightarrow L $ is \emph{homogeneous} if it is invariant under continuous mappings of the form  $\gamma (x)=x\wedge c$ and $\gamma (x)=x\vee c$, for every $c\in L$. In particular, if $F$ is homogeneous, then it is idempotent.  Moreover, as in the case of discrete Sugeno integrals (see \cite{CouMar}), we have the following result.

\begin{proposition}\label{Proposition:homo} Every Sugeno integral (in particular, every term functional)  is homogeneous and nondecreasing.
\end{proposition}

\begin{proof} As observed, every polynomial functional and, in particular, every Sugeno integral is nondecreasing. 
By Theorem~\ref{Prop:2}, if $F \colon  L^{A}\rightarrow L $ is a Sugeno integral, then for every $f\in L^A$ we have
$$F(f)= \bigvee_{X\in\mathcal A} F(I_X) \wedge \bigwedge_{x\in X} f(x)$$
where $F(I_\emptyset)=0$ and $F(I_A)=1$. Thus, if $\gamma \colon L\to L$ is a continuous mapping of the form  $\gamma (x)=x\wedge c$ or $\gamma (x)=x\vee c$, for some $c\in L$, it follows from complete distributivity that 
\begin{eqnarray*}
\gamma \circ F(f)&=& \gamma (\bigvee_{X\in\mathcal A} F(I_X) \wedge \bigwedge_{x\in X} f(x)) \\
&=&
\bigvee_{X\in\mathcal A}  \gamma (F(I_X) \wedge \bigwedge_{x\in X} f(x))= 
\bigvee_{X\in\mathcal A} F(I_X) \wedge \bigwedge_{x\in X}  \gamma (f(x)).
\end{eqnarray*}
In other words, $F$ is homogeneous. 
\end{proof}

Observe that Proposition~\ref{Proposition:homo} does not hold if homogeneity is replaced by invariance. To see this, let $L$ be the real interval $[0,1]$ and 
fix $c\in (0,1)$. Consider the Sugeno integral 
$$F(x_1, x_2)=\med(x_1,c,x_2)=(x_1\wedge c)\vee (x_1\wedge x_2)\vee (c\wedge x_2),$$  
and let $\gamma (x)=x^2$. Clearly, $\gamma$ is continuous but
$$\gamma \circ F(x_1, x_2)=\med(x_1^2,c^2,x_2^2)\neq \med(x_1^2,c,x_2^2)=F(\gamma(x_1),\gamma(x_2)).$$
However, we have the following result. 
 
\begin{proposition}\label{Proposition:inv} Every term functional is invariant and nondecreasing.
\end{proposition}

\begin{proof}
By Corollary~\ref{cor:1},  every term functional $F \colon  L^{A}\rightarrow L $ can be represented by
$$F(f)=\bigvee_{X\in\mathcal A} \bigwedge_{x\in X} f(x)=\bigwedge_{X\in\mathcal B} \bigvee_{x\in X} f(x)=P^\mathcal{B}(f),$$
where $\mathcal A=\{X\subseteq A\colon F(I_X)=1\}$ and $\mathcal B=\{X\subseteq A\colon F(I_{A\setminus X})=0\}$.
Thus, if $\gamma \colon L\to L$ is a continuous mapping, then 
$$
\gamma \circ F(f)= \gamma (\bigvee_{X\in\mathcal A} \bigwedge_{x\in X} f(x))= \bigvee_{X\in\mathcal A} \gamma (\bigwedge_{x\in X} f(x))
= \bigvee_{X\in\mathcal A}\bigwedge_{x\in X} \gamma (f(x)).
$$
In other words, $F$ is invariant. 
 \end{proof}

\subsection{Term functionals as invariant functionals on completely distributive lattices}

As observed, every term functional is invariant and nondecreasing. The following result shows that these two conditions are in fact sufficient to axiomatize this subclass of Sugeno integrals on completely distributive lattices.

\begin{theorem}\label{Prop:3.2a}
Let $L$ be a completely distributive lattice, $A$ an arbitrary
nonempty set, and let $F \colon  L^{A}\rightarrow L $ be a functional such that, for every $X\subseteq A$, $F(I_X)\in \{0,1\}$. 
Then $F$ is a term functional if and only if
it is nondecreasing and invariant.
\end{theorem}

\begin{proof}
By Proposition~\ref{Proposition:inv}, every term functional is both invariant and nondecreasing.

To show that the converse also holds, let $F \colon  L^{A}\rightarrow L $ be a nondecreasing invariant functional 
such that, for every $X\subseteq A$, $F(I_X)\in \{0,1\}$. We show that
 $P_F(f)\leqslant F(f)\leqslant P^F(f) $, for every $f\in L^A$, and thus that $F$ is a polynomial functional by Theorem \ref{Prop:2}.

 For each $X\subseteq {A}$, consider the continuous mapping $\gamma_X $ given by
  $\gamma_X (x) = x \wedge (\bigwedge f(X))$. Since $F$ is invariant, we have
  $F(I_X)\wedge \bigwedge f(X)=\gamma_X\circ F(I_X)=F(\gamma_X\circ I_X)$,
  and since $F$ is nondecreasing, $F(I_X)\wedge \bigwedge f(X)\leqslant F(f)$. Thus $P_F(f)\leqslant F(f)$.

 Now, for each $Y\subseteq {A}$, consider the continuous mapping $\gamma^Y $ given by
    $\gamma ^Y (x) = x \vee (\bigvee f(Y))$. Then,
    $F(I_{A\setminus Y})\vee \bigvee f(Y)= \gamma ^Y\circ F(I_{A\setminus Y})= F(\gamma ^Y\circ I_{A\setminus Y}) $,
    and since $F$ is nondecreasing, $F(I_{A\setminus Y})\vee \bigvee f(Y)\geqslant F( f)$. Thus $P^F(f)\geqslant F(f)$.
  
  Since $F(I_X)\in \{0,1\}$, for every $X\subseteq A$, we have that $F$ is a term functional, and the proof is now complete.  
\end{proof}

\begin{remark}
A variant of Theorem \ref{Prop:3.2a} was first obtained by Ovchinnikov and Dukhovny \cite{OvDuck02} in the particular case when $L$ is a complete dense chain.
\end{remark}

In the case  when $L$ is a complete chain, Theorem \ref{Prop:3.2a} can be strengthened since nondecreasing monotonicity becomes redundant as the following result asserts.

\begin{proposition}\label{Prop:3.3}
Let  $L \neq \{ 0, 1\}$ be a complete chain and let $A$ be a nonempty
set. Then every invariant functional $F \colon L^{A}\rightarrow L $
 is  nondecreasing.
\end{proposition}

\begin{proof} Let  $F \colon L^{A}\rightarrow L $ be an invariant functional, and for the sake of contradiction, suppose that there exist $ f, g \in L^{A}$ such that $f\leq g$ but $F(f)>F(g)$.

{\it Case 1:} There exists $t\in L$ such that $F(f)>t>F(g)$.

Let $\gamma ^t\colon L\to L$ and $\gamma _t\colon L\to L$ be the continuous mappings given by
$$\gamma ^t(x)=x\vee t \quad \mbox{ and } \quad \gamma _t(x)=x\wedge t, \quad \mbox{resp.},$$
and set $h =(\gamma^t\circ f) \wedge g$. Then we have
\begin{eqnarray*}
\gamma^t\circ h &=& (\gamma^t\circ f) \wedge (\gamma^t\circ g ) = \gamma^t\circ f,\\
\gamma_t\circ h &=& t \wedge (\gamma_t\circ g ) = \gamma_t\circ g.
\end{eqnarray*}
Since  $F$ is invariant, and $\gamma ^t$ and $\gamma _t$ are continuous,
$ F(h) \vee t = F(f) \vee t = F(f)$ and hence $ F(h) = F(f)$ because $L$ is a chain and $ F(f)> t$.
On the other hand, we also have $ F(h) \wedge t = F(g) \wedge t = F(g)$ and hence $ F(h) = F(g)$, which yields the desired contradiction.

{\it Case 2:} There is no $t\in L$ such that $F(f)>t>F(g)$.

Let $f'=\med(F(f),f, F(g))$ and $g'=\med(F(f),g, F(g))$. By invariance, it follows that $F(f')=F(f)>F(g)=F(g')$ and since $f<g$, we have $f'<g'$. Thus, we may assume that $f=f'$ and $g=g'$.
In other words, there are $ K'\subsetneq K\subseteq A$  such that
\[
f(a) =
\begin{cases}
F(g) & \text{if $a\in K$} \\
F(f), & \text{otherwise,}
\end{cases}
\]
\[
g(a) =
\begin{cases}
F(g) & \text{if $a\in K'$} \\
F(f), & \text{otherwise.}
\end{cases}
\]

Observe that, for every $h_g\in L^A$ such that $h_g(a) =g(a)$, if $a\in K'$, and $h_g(a)\geqslant g(a)$, otherwise, we must have $F(h_g)= F(g)$ because $F$ is invariant under $x\wedge F(f)$.

Similarly, for every $h_f\in L^A$ such that $h_f(a) =f(a)$, if $a\in K$, and $h_f(a)\geqslant f(a)$, otherwise, we have $F(h_f)\geqslant F(f)$.
In fact, for every such $h_f$, we must have $F(h_f)= F(f)$. Indeed, if there is an $b\in A\setminus K$ such that $h_f(b)>f(b)$ and $F(h_f)> F(f)$, then  by taking $t=F(f)$ and $h_g\in L^A$ such that $h_g(a) =g(a)$, if $a\in K'$, and $h_g(a)=h_f(b)$, otherwise, we would be in Case 1 since $h_f<h_g$ and $F(h_f)>t> F(h_g)$, and thus obtain the same contradiction.

Suppose that $F(f)\neq 1$, i.e., $F(f)< 1$.  Define $\gamma\colon L\to L$ by $\gamma (x)=x$, if $x\leqslant F(g)$, and $\gamma(x)=1$, otherwise. Since there is no $t\in L$ such that $F(f)>t>F(g)$, we have that $\gamma $ is continuous. However,
$$F(\gamma\circ f)=F(f)<1=\gamma\circ F(f)$$
which contradicts the fact that $F$ is invariant.

If $F(f)= 1$, then $F(g)>0$ and for $h_f^0\in L^A$ such that $h_f^0(a) =0$, if $a\in K$, and $h_f^0(a)= f(a)$, otherwise, we have $F(h_f^0)\geqslant F(f)=1$.
Also, observe that for every $h_g\in L^A$ such that $h_g(a) \leqslant g(a)$, if $a\in K'$, and $h_g(a)=g(a)$, otherwise, we must have $F(h_g)\leqslant F(g)$ because $F$ is invariant under $x\vee F(g)$.

Now, if there is such $h_g$ for which $F(h_g)<F(g)<1=F(f)$, then we can proceed as in Case 1 with $t=F(g)$ since $h_g>h_f^0$ and $F(h_f^0)=F(f)$. Otherwise, taking $h_g^0\in L^A$ such that
 $h_g^0(a) =0$, if $a\in K'$, and $h_f^0(a)= f(a)$, otherwise, and defining $\gamma\colon L\to L$ by $\gamma (x)=0$, if $x\leqslant F(g)$, and $\gamma(x)=x$, otherwise, we have
 $$F(\gamma\circ g)=F(g)>0=\gamma\circ F(g)$$
 which yields the desired contradiction.
\end{proof}

From Theorem \ref{Prop:3.2a} and Proposition \ref{Prop:3.3}, we obtain the following characterization of term functionals over complete chains  $L\neq \{ 0, 1\}$.
 
\begin{corollary}
Let $L\neq \{ 0, 1\}$ be a complete chain, $A$ an arbitrary
nonempty set, and let $F \colon  L^{A}\rightarrow L $ be a functional such that, for every $X\subseteq A$, $F(I_X)\in \{0,1\}$. 
Then $F$ is a term functional if and only if
it is invariant.
\end{corollary}

Observe that Proposition \ref{Prop:3.3} does not hold when $L=\{ 0, 1\}$. Indeed, in this case, every idempotent functional is invariant, but one can easily find examples of idempotent Boolean functions which are not nondecreasing.
Also, invariance of $F : L^{A}\rightarrow L $ cannot be replaced by the weaker condition of homogeneity.
To illustrate, let $L=\{0,a,1\}$ be a chain, $A=\{1,2,3\}$, and consider the symmetric functional $F : L^{A}\rightarrow L $
given by
\begin{eqnarray*}
F(f(1),f(2),f(3)) =
\begin{cases}
x & \text{if $f(1)=f(2)=f(3)=x$,} \\
a & \text{if }f(i)=f(j)=0,\, f(k)=a,1, \text{ and}\\
   & \{i,j,k\}=A, \\
f(i) & \text{if }f(i)\in\{0,a\},\, f(j),f(k)\in \{a,1\},\text{ and}\\
      & \{i,j,k\}=A.
\end{cases}
\end{eqnarray*}
Clearly, $F$ is homogeneous but it is not nondecreasing.

\subsection{Characterizations of Sugeno integrals and polynomial functionals on completely distributive lattices}

As we have seen, not every Sugeno integral (and thus not every polynomial function) is invariant. However, Proposition~\ref{Proposition:homo} asserts that Sugeno integrals fulfill the weaker invariance property of homogeneity.
As the following result shows, under nondecreasing monotonicity, this homogeneity condition suffices to characterize Sugeno integrals on completely distributive lattices.

\begin{theorem}\label{Prop:3.2}
Let $L$ be a completely distributive lattice, $A$ an arbitrary
nonempty set, and let $F \colon  L^{A}\rightarrow L $ be a functional. 
Then $F$ is a Sugeno integral if and only if
it is nondecreasing and homogeneous.
\end{theorem}

\begin{proof}
By Proposition~\ref{Proposition:homo}, the conditions are necessary.

To show that the converse also holds, let $F \colon  L^{A}\rightarrow L $ be a nondecreasing homogeneous functional.
 %We follow exactly the same steps as in the proof of Theorem~\ref{Prop:3.2a}, and show 
 First, we show  that
 $P_F(f)\leqslant F(f)\leqslant P^F(f) $, for every $f\in L^A$. 

 For each $X\subseteq {A}$, consider the continuous mapping $\gamma_X $ given by
  $\gamma_X (x) = x \wedge (\bigwedge f(X))$. Since $F$ is homogeneous, we have
  $F(I_X)\wedge \bigwedge f(X)=\gamma_X\circ F(I_X)=F(\gamma_X\circ I_X)$,
  and since $F$ is nondecreasing, $F(I_X)\wedge \bigwedge f(X)\leqslant F(f)$. Thus $P_F(f)\leqslant F(f)$.

 Now, for each $Y\subseteq {A}$, consider the continuous mapping $\gamma^Y $ given by
    $\gamma ^Y (x) = x \vee (\bigvee f(Y))$. Then,
    $F(I_{A\setminus Y})\vee \bigvee f(Y)= \gamma ^Y\circ F(I_{A\setminus Y})= F(\gamma ^Y\circ I_{A\setminus Y}) $,
    and since $F$ is nondecreasing, $F(I_{A\setminus Y})\vee \bigvee f(Y)\geqslant F( f)$. Thus $P^F(f)\geqslant F(f)$.
   
   By Theorem~\ref{Prop:2}, we have that $P_F(f)= P^F(f) $, for every $f\in L^A$ and hence $F=P_F= P^F $.
 Moreover, from homogeneity it follows that $F$ is idempotent. Thus $F$ is a Sugeno integral, and 
  the proof of the theorem is now complete.  
\end{proof}

The superclass of lattice polynomial functionals can be similarly axiomatized but in terms of a weaker form of homogeneity.
We say that a functional  $F\colon L^A\to L$ is \emph{range-homogeneous} if, for every $f\in F^A$,
$$F(\gamma \circ f)=\gamma \circ F(f),$$
whenever $\gamma \colon L\to L$ is a continuous mapping of the form  
$\gamma (x)=x\wedge c$ or $\gamma (x)=x\vee c$, and where $c\in \overline{\mathcal{R}}_F $. 
Note that every lattice polynomial functional $F\colon L^A\to L$ is nondecreasing and range-homogeneous.
In particular, for every $c\in  \overline{\mathcal{R}}_F$, $F(\vect{c})=c$.

Based on this weaker notion of homogeneity, we have the following characterization of lattice polynomial functionals.

\begin{theorem}\label{Prop:3.22}
Let $L$ be a completely distributive lattice, $A$ an arbitrary
nonempty set, and let $F \colon  L^{A}\rightarrow L $ be a functional. Then $F$ is a lattice polynomial functional
if and only if it is nondecreasing and range-homogeneous.
\end{theorem}

\begin{proof} As observed, the conditions are necessary.
 The proof of the sufficiency follows the same steps as in the proof of Theorem \ref{Prop:3.2},
by using Theorem \ref{Prop:2} after showing that $P_F(f)\leqslant F(f)\leqslant P^F(f) $.

For every $X, Y\subseteq A$, define 
\begin{eqnarray*}
\gamma_X(x)&=&\langle x \wedge (\bigwedge f(X))\rangle_F~= \langle x \rangle_F \wedge \langle(\bigwedge f(X))\rangle_F,\quad \mbox{ and }\\ \gamma^Y(x)&=&\langle x  \vee (\bigvee f(Y))\rangle_F~=~ \langle x \rangle_F\vee \langle(\bigvee f(Y))\rangle_F,
\end{eqnarray*}
where the operator $\langle \cdot \rangle_F\colon L\to L$ is given by
$$\langle x \rangle_F=\med(\bigwedge \overline{\mathcal{R}}_F,x,\bigvee  \overline{\mathcal{R}}_F)=
( \bigwedge \overline{\mathcal{R}}_F\vee x)\wedge \bigvee  \overline{\mathcal{R}}_F=
\bigwedge \overline{\mathcal{R}}_F\vee (x\wedge \bigvee  \overline{\mathcal{R}}_F).$$

By range-homogeneity, $F(\langle f\rangle_F)=\langle F(f)\rangle_F=F(f)$, for all $f\in L^A$, and 
\begin{eqnarray*}
 F(\gamma_X\circ I_X)&=&\gamma_X\circ F( I_X)= \langle F( I_X) \rangle_F \wedge \langle \bigwedge f(X)\rangle_F,\quad \mbox{ and }\\
F(\gamma ^Y\circ I_{A\setminus Y})&=& \gamma^Y\circ F(I_{A\setminus Y})=\langle F(I_{A\setminus Y}) \rangle_F\vee \langle \bigvee f(Y)\rangle_F,
\end{eqnarray*}
 for all $X, Y\subseteq A$. Moreover, by nondecreasing monotonicity, it follows that for every $X, Y\subseteq A$ and $f\in L^A$
\begin{eqnarray*}
F(f)= F(\langle f\rangle_F) &\geqslant& F(\gamma_X\circ I_X) =  \langle F( I_X) \rangle_F \wedge \langle \bigwedge f(X)\rangle_F,\quad \mbox{ and }\\
F(f)= F(\langle f\rangle_F) &\leqslant& F(\gamma ^Y\circ I_{A\setminus Y}) =\langle F(I_{A\setminus Y}) \rangle_F\vee \langle \bigvee f(Y)\rangle_F.
\end{eqnarray*}
Hence, for every $f\in L^A$, we have
\begin{eqnarray*}
P_F(f)&=&\langle P_F(f)\rangle_F
= \bigvee_{X\subseteq A} \langle F( I_X) \rangle_F \wedge \langle \bigwedge f(X)\rangle_F\leqslant F(f) \\
&\leqslant&  \bigvee_{Y\subseteq A} \langle F(I_{A\setminus Y}) \rangle_F\vee \langle \bigvee f(Y)\rangle_F
= \langle P^F(f) \rangle_F = P^F(f).
 \end{eqnarray*}
 Since $P_F= P^F $ (Theorem~\ref{Prop:2}), $F$ is a lattice polynomial functional, and the proof of the theorem is now complete.
 \end{proof}

\begin{remark}
In the particular case when $A$ is finite, Theorems~\ref{Prop:3.2} and \ref{Prop:3.22} reduce to Theorem 24 $(iv)$ in \cite{CouMar} and Main Theorem $(iv)$ in \cite{CouMar1}, respectively, since in this case it is sufficient to require that $L$ is bounded and distributive. 
\end{remark}

\section*{Acknowledgments}
We would like to thank Radko Mesiar for helpful suggestions concerning Subsection 4.1, in particular, for providing the example showing that not every Sugeno integral is invariant.

\end{document}